\documentclass[12pt]{amsart}
\usepackage[mathscr]{eucal}

\usepackage{amsfonts}
\usepackage{amsmath}
\usepackage{amsthm}
\usepackage{amssymb}
\usepackage{latexsym}
\usepackage{epsfig}
\usepackage{graphicx}
\usepackage{hhline}

\usepackage{pstricks}

\newtheorem{theorem}{Theorem}
\newtheorem{lemma}[theorem]{Lemma}
\newtheorem{remark}[theorem]{Remark}

\newtheorem{definition}{Definition}

\newtheorem{proposition}[theorem]{Proposition}

\newcommand{\eps}{\varepsilon}
\newcommand{\R}{{\mathbb{R}}}
\newcommand{\Q}{{\mathbb{Q}}}
\newcommand{\Z}{{\mathbb{Z}}}
\newcommand{\N}{{\mathbb{N}}}
\newcommand{\T}{{\mathbb{T}}}

\newcommand{\Prob}{\mathbb{P}}

\newcommand{\qq}{\mathcal{Q}}
\newcommand{\qqe}{\mathcal{Q}_\eps}
\newcommand{\qqexa}{\mathcal{Q}_\eps (x,\alpha)}

\newcommand{\nn}{\mathcal{N}}
\newcommand{\rect}{\mathcal{R}}

\newcommand{\aaa}{\mathcal{A}}

\newcommand{\ua}{\left(\begin{array}{ll} 1 & \alpha \\ 0 & 1 \end{array}\right)}

\newcommand{\dnet}{\left(\begin{array}{ll} \eps T^{-1} & 0\\ 0 & \eps^{-1} T \end{array}\right)}

\newcommand{\de}{\left(\begin{array}{ll} \eps & 0\\ 0 & \eps^{-1} \end{array}\right)}

\newcommand{\ben}{\begin{enumerate}}
\newcommand{\een}{\end{enumerate}}
\newcommand{\beq}{\begin{equation}}
\newcommand{\eeq}{\end{equation}}

\newcommand{\vo}{v_{\rm out}}
\newcommand{\vin}{v_{\rm in}}
\newcommand{\tht}{\hat{T}}
\newcommand{\thtae}{\hat{T}_{\alpha, \eps}}

\newcommand {\vphi}  {\varphi}
\newcommand {\Om}   {\Omega}
\newcommand{\xe}{\xi_\eps}

\newcommand{\bfy}{ {\bf y} }

\begin{document}

\title{Perfect Retroreflectors and Billiard Dynamics}

\author[P. Bachurin]{Pavel Bachurin}
\address{Department of Mathematics, SUNY Stony Brook, USA}
\email{bachurin@math.toronto.edu}

\author[K. Khanin]{Konstantin Khanin}
\address{Department of Mathematics, University of Toronto, Canada}
\email{khanin@math.toronto.edu}
\thanks{KK is supported by an NSERC Discovery grant}

\author[J. Marklof]{Jens Marklof}
\address{Department of Mathematics, University of Bristol, UK}
\email{j.marklof@bristol.ac.uk}
\thanks{JM is supported by a Royal Society Wolfson Research Merit Award}

\author[A. Plakhov]{Alexander Plakhov}
\address{Department of Mathematics, University of Aveiro, Portugal} 
\email{plakhov@ua.pt}
\thanks{AP is supported by {\it Centre for Research on Optimization
and Control} (CEOC) from the ''{\it Funda\c{c}\~{a}o para a
Ci\^{e}ncia e a Tecnologia}'' (FCT), cofinanced by the 
European Community Fund FEDER/POCTI, and by the FCT 
research project PTDC/MAT/72840/2006}



\begin{abstract}

We construct semi-infinite billiard domains which reverse the direction of most incoming particles.
We prove that almost all particles will leave the open billiard domain after a finite number of reflections.
Moreover, with high probability the exit velocity is exactly opposite to the entrance velocity, and the 
particle's exit point is arbitrarily close to its initial position. The method is based on asymptotic analysis 
of statistics of entrance times to a small interval for irrational circle rotations. The rescaled entrance 
times have a limiting distribution in a limit when the number of iterates tends to infinity and the length 
of the interval vanishes. The proof of the main results follows from the study of related limiting distributions 
and their regularity properties.

\end{abstract}


\maketitle

\renewcommand{\theequation}{\arabic{section}.\arabic{equation}}

\section{Introduction}

The present paper is motivated by the problem of constructing open billiard domains with exact velocity reversal (EVR), which means that the 
velocity of every incoming particle is reversed when the particle eventually leaves the domain. This problem arises in the construction 
of perfect retroreflectors---optical devices that exactly reverse the direction of an incident beam of light and preserve the original image. A 
well-known example of a perfect retroreflector is the Eaton lens \cite{E}, \cite{T} which is a spherically symmetric lens that, unlike our model, 
also reverses the original image. A second application lies in the search for domains that maximize the pressure of a flow of particles \cite{PG}: 
for a particle of mass $m>0,$ which moves towards a wall with velocity $\bar{v},$  the impulse transmitted to the wall at the moment of reflection 
is equal to $2m |\bar{v}_n|,$ where $\bar{v}_n$ is the normal component of $\bar{v}.$ It is maximized when  $\bar{v}=\bar{v}_n,$ i.e. when 
the direction of the particle is reversed.  

We construct a family of domains $D_\epsilon$, for which EVR holds up to a set of initial condition whose measure tends to zero in the 
limit $\eps\to 0$.

The domain $D_\epsilon$ is the semi-infinite tube $[0,\infty)\times [0,1]$ with vertical barriers of height $\epsilon/2$ at the points 
$(n, 0)$ and $(n,1),$ $n\in\N$ as illustrated in Fig. 1. Inside the domain the particle moves with the constant  speed and elastic reflections 
from the boundary. Since the kinetic energy of the particle is preserved, we can assume that the speed of the particle is equal to one.

The motion of the particle is determined by the point $y_{\rm in}\in [0,1] ,$ where it enters the tube and the initial velocity 
$\vin=(\cos(\pi\vphi),\sin(\pi\vphi))$ at this point. The measure $\Prob$ on the initial conditions $( y_{\rm in}, \vphi)$ considered below 
is a Borel probability measure absolutely continuous with respect to the Lebesgue measure on $\Om= [0,1]\times [-1/2,1/2].$


%
%

\bigskip

\begin{figure} [h]
	\centering
	\includegraphics [width=4in] {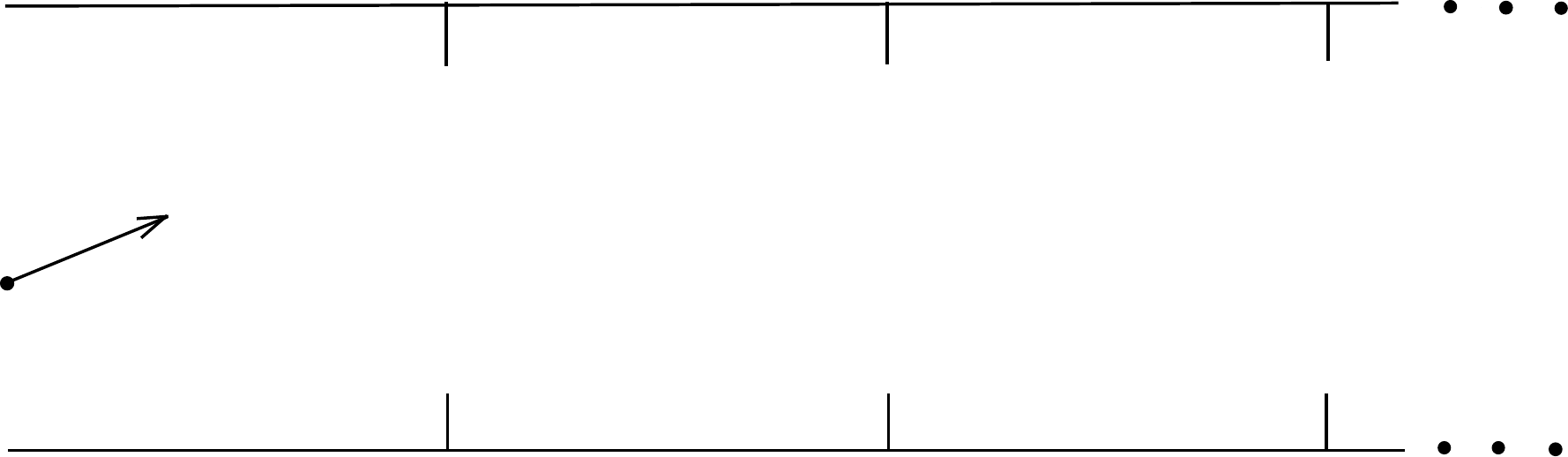}
	\caption{The Model}
	\label{Fig1}
\end{figure}

\begin{theorem}
    For every $\eps\in (0,1)$ there exists a set $\Om(\eps)\subset\Om$ of  full Lebesgue measure, 
    such that for every $( y_{\rm in}, \vphi)\in\Om(\eps),$ the particle eventually leaves the tube.
    \label{asexits}
\end{theorem}

The position and the velocity with which it leaves the tube are denoted by $(y_{\rm out}, \vo).$ By
Theorem \ref{asexits}, for every $\eps\in (0,1)$ the functions $y_{\rm out}=y_{\rm out}(y_{\rm in},v_{\rm in})$ and 
 $\vo=\vo( y_{\rm in}, \vin)$ are defined $\Prob$-almost everywhere.  
\begin{theorem} For any $\delta>0,$
	\beq
		\Prob\{(y_{\rm in}, \vin)\ :\ \vo=-\vin,\ |y_{\rm out}-y_{\rm in}|<\delta\}\to 1\quad {\rm as}\ \eps\to 0
	\label{goal}
	\eeq
	\label{maintheorem}
\end{theorem}
Theorem \ref{maintheorem} follows from the results on the existence of certain limiting distributions for 
the exit statistics of the billiard particle as $\eps\to 0.$ Below we formulate these results as Theorem \ref{thm1}
and Theorem \ref{thm2}. In the last section of the paper we show how they imply Theorem \ref{maintheorem}.

Let $\qqe=\qqe(y_{\rm in},\vin)$ be the number of reflections from the vertical walls  before the particle  
leaves the tube. Let $T_\eps=T_\eps(y_{\rm in},\vin)$ be the time that particle spends inside the tube.
By Theorem \ref{asexits}, both $\qqe$ and $T_\eps$ are finite $\Prob-$a.e.

Consider also a bi-infinite tubular domain similar to the one described above. It consists of two horizontal lines 
at the unit distance from each other and a one-periodic configuration of vertical walls of height $\eps/2.$

Let $x$ be the  horizontal coordinate and assume that the particle starts at $x=0.$  Let $\xe^0=0$ and 
$\xe^k\in\Z$ be $x$-coordinate of the particle at the moment of $k$'th reflection from a vertical wall. Since 
the tube is now bi-infinite, $\{\xe^k\}$ is a discrete time process on $\Z,$ defined for any $k\in\N.$ We also 
define a continuous version of this process: $\{\xe(t)\}$ is the projection of the trajectory of a billiard particle 
in the bi-infinite tube to the $x$-axis normalized to have constant speed $1/\eps.$

\begin{theorem}
    \begin{enumerate}
	 \item The process $\{\eps \xe^k\}$ converges in distribution (w.r.t. $\Prob$) to a stochastic process $\{\xi^k\}$ as $\eps\to 0.$    
          \item There exists a limiting probability distribution function $G: \N\to [0,1]$ such that for every $k\in\N,$ $\Prob \{\qqe(y_{\rm in}, \vin) = k \}\to G(k)$ as $\eps\to 0$
    \end{enumerate}
    \label{thm1}
\end{theorem}

The second part of Theorem \ref{thm1} says that for the limiting stochastic process $\{\xi^k\}$, with probability one there exists 
$k\in\N,$ such that  $\xi^k <0.$ Similar results are true for the continuous process $\{ \xe(s)\}$ as well:

\begin{theorem}
	\begin{enumerate}
		\item The process $\{\eps \xe(s)\}$ converges in distribution w.r.t. $\Prob$ to a stochastic process $\xi(s)$ as $\eps\to 0.$ 		
		\item There exists a limiting probability distribution function $H:\R_{\geq 0}\to [0,1],$  such that for every $t\geq 0,$ $\Prob\{\eps T_\eps(y_{\rm in},\vin) < t\}\to H(t)$ as $\eps\to 0.$
	\end{enumerate}	
	\label{thm2}
\end{theorem}


\section {Reduction to Circle Rotations and Point-Wise Exits}
\label{point-wise}

We first reformulate the problem in terms of circle rotations.

Let us identify $[0,1)$ with $S^1=\R/\Z.$ For $\alpha \in \R,$ let $R_\alpha:S^1\to
S^1$ be the circle rotation by angle $\alpha:$
$$R_\alpha x= x+\alpha \mod 1.$$

Always assume that $\alpha\in\R\setminus\Q.$

Let $I_\eps=[-\eps/2,\eps/2]\subset S^1.$ 

We define several sequences measuring the return times to the interval $I_\eps,$ which 
will be used throughout the proofs. The hitting times $m_{\eps}^k=m_{\eps}^k(x,\alpha),$ 
$k=0,1,2,\ldots$ are defined for $x\in S^1$ by:
$$
    m_{\eps}^0=0,\ m_{\eps}^k(x)=\min\{l>m_{\eps}^{k-1}:\ R_\alpha^l x\in I_\eps\}
$$

The sequence $n_\eps^k=n_{\eps}^k(x,\alpha),$ $k=1,2,\ldots$  of  {\it relative} return times to the interval $I_\eps$
is defined for $x\in S^1$ by:
$$
    n_{\eps}^k=m_{\eps}^k(x)-m_\eps^{k-1}(x)
 $$

We shall also use the sequence $\{\xi_{\eps}^k\}$  defined in the introduction as the sequence of the horizontal coordinates 
of points of the reflection from the vertical walls.

Note that if $x=y_{\rm in},$ and $\alpha=\tan(\pi\vphi),$ then $n^i_{\eps}(x)$ is the distance between
horizontal coordinate of the place of the $(i-1)$'st and the $i$'th reflections of the particle from vertical walls.
Therefore,
$$
    		\xi_\eps^k=n_\eps^1-n_\eps^2+\ldots+(-1)^{k+1}n_\eps^k,\\ \\
$$
$$
    		\qqe=\qqexa=\min\{j\in \N\ :\ n_{\eps}^1(x)-n_{\eps}^2(x)+\ldots + (-1)^{j+1} n_{\eps}^j(x)  \leq 0\} -1
$$

Let $\bar{n}_{\eps}^k=(n_{\eps}^1,\ldots, n_{\eps}^k)^T,$ $\bar{m}_{\eps}^k=(m_{\eps}^1,\ldots, m_{\eps}^k)^T$ and 
$\bar{\xi}_{\eps}^k=(\xi_{\eps}^1,\ldots, \xi_{\eps}^k)^T,$  then 
\beq
	 \bar{\xi}_{\eps}^k= {\bf A} \bar{n}_{\eps}^k,\ {\rm and}\ \bar{m}_{\eps}^k={\bf B} \bar{n}_{\eps}^k,
\label{relation}
\eeq	 

where ${\bf A}$ and ${\bf B}$ are two $k\times k$ matrices with

$${\bf A}_{i,j}=
	\left\{\begin{array}{l} 0,\ {\rm if}\ i<j, \\ (-1)^{j+1},\ {\rm if}\ i\geq j \end{array}\right. ,\quad {\rm and}\quad
     {\bf B}_{i,j}=	
	\left\{\begin{array}{l} 0,\ {\rm if}\ i<j, \\  1,\ {\rm if}\ i\geq j \end{array}\right. 
$$

The probability measure $\Prob$ on the initial conditions $(y_{\rm in},\vphi_{in})\in [0,1]\times [-1/2,1/2]$ for the billiard particle induces
a probability measure on the initial conditions $(x,\alpha)\in [0,1)\times [0,1)\simeq \T^2$ for the circle rotation $R_\alpha,$ which is absolutely 
continuous w.r.t. to the Lebesgue measure on $\T^2$ and which will be also denoted by $\Prob.$

We now prove Theorem \ref{asexits}. 

Let $\tht_{\alpha,\eps}:I_{\eps}\to I_{\eps}$ be the map induced on $I_\eps$ by the circle
rotation $R_\alpha:$
$$\tht_{\alpha,\eps}(x)=R_{\alpha}^{m_\eps^1}(x)$$
 

\begin{proposition} For every $\eps\in (0,1)$ there exists a set of full Lebesgue measure
    $\Lambda(\eps)\subset S^1,$ such that for every $\alpha\in\Lambda,$ the map 
    $\tht_{\alpha,\eps}$ is weakly mixing.
\label{iet}
\end{proposition}

\begin{proof}

The proof of Proposition will follow from a combination of results of \cite{B} and \cite{BN}.

For every $\eps>0$ there exists a full Lebesgue measure set $\Lambda '(\eps)\subset S^1,$  such that
for every $\alpha\in\Lambda '(\eps),$ the corresponding map  $\tht_{\alpha,\eps}$ is an interval exchange transformation 
of three intervals of combinatorial type $(3,2,1).$

Recall Property P introduced by Boshernitzan in \cite{B}:

\begin{definition}
 A set $\aaa\subset\N$ is called essential, if for any integer $l\geq 2$ there exists $\lambda >1,$ such that 
the system 
\[\left\{\begin{array} {ll}
	n_{i+1}> 2 n_i, & {\rm for}\ 1\leq i\leq l-1,\\
	n_l< \lambda n_1, & \\
	n_i\in\aaa, & {\rm for}\ 1\leq i\leq l
\end{array}
\right.
\]
has an infinite number of solutions $(n_1, n_2,\ldots n_l).$
\end{definition}
 
Let $m_n(\thtae)$ be the length of the smallest interval of continuity of $\thtae^n.$

\begin{definition} An interval exchange map $\tht_{\alpha, \eps}$ has property P, if for some $\delta>0$
the set $$\aaa(\alpha,\eps, \delta)=\{n\in\N | m_n(\thtae) > \frac{\delta}{n} \}$$ is essential.
\end{definition}

\begin{proposition} (\cite{B}, Theorem 9.4 (a)) For every $\eps\in (0,1)$ there exist a full Lebesgue measure
set $\Lambda(\eps)\subset S^1,$ such that for every $\alpha\in\Lambda,$ the map $\tht_{\alpha,\eps}$  has property P.
\end{proposition}

By Theorem 5.3 of \cite{BN}, property P implies weak-mixing for an interval exchange of three
intervals with combinatorics $(3,2,1)$ (and more generally, for any combinatorics of a so-called $W-$type,
see \cite{BN}).

This implies Proposition \ref{iet}.
\end{proof}

The next two statements are well-known. We include their proofs to keep the exposition self-contained.

\begin{lemma} For every $\eps\in (0,1)$ there exists a set of full Lebesgue measure
    $\Lambda(\eps)\subset S^1,$ such that for every $\alpha\in\Lambda,$ the map $\tht_{\alpha,\eps}^2$ is ergodic.
\label{square}
\end{lemma}

\begin{proof}
Assume now that $\tht_{\alpha,\eps}^2$ is not ergodic. Then there exists a bounded $f\ne const,$
such that $\tht_{\alpha,\eps}^2 f=f,$ and therefore $$\tht_{\alpha,\eps}(f+\tht_{\alpha,\eps} f)=
\tht_{\alpha,\eps} f + f$$

Since $\tht_{\alpha,\eps}$ is ergodic, this implies that $\tht_{\alpha,\eps} f + f = C,$ or
$f-C/2=-(\tht_{\alpha,\eps} f - C/2).$ If $g=f-C/2,$ then $g$ is not identically zero,
and $\tht_{\alpha,\eps} g=-g.$ Therefore $\lambda=-1$ is an eigenvalue of
$\tht_{\alpha,\eps}$ and so $\tht_{\alpha,\eps}$ is not weakly mixing.
\end{proof}

\begin{proposition} Let $T$ be an ergodic transformation on $(X,\mu),$ $\mu(X)=1,$ and
let $f\in L^1(X,\mu),$ $\int f d\mu=0$ and
$S_n(f,x)=f(x)+f(Tx)+\ldots+f(T^{n-1}x)$ be its Birkhoff sums. Then either
$S_n(f,x)$ is unbounded from below for almost every $x\in X,$ or $f$ is a
co-boundary, i.e. there exists a measurable $g(x)$, such that $f(x)=g(x)-g(Tx).$
\label{abs}
\end{proposition}

\begin{proof} Since $T$ is ergodic, the set of points $x$ for which $S_n(f,x)$ is bounded from 
below has measure either equal to zero or one. In the first case, Proposition is proved, so  
assume that it has measure one.  Then the function  $g(x)=\inf\limits_{n\geq 1} S_n(f,x)$ is 
finite almost everywhere.

We have $g(Tx)+f(x)=\inf\limits_{n\geq 2} S_n(f,x),$ so $h(x):=g(Tx)-g(x)+f(x)\geq 0.$

It is enough to show that $\int h d\mu=0.$  If $g(x)\in L^1(X,\mu),$ then $\int h d\mu=0$ by the definition of $h(x)$ above. 
If not, then by Birkhoff ergodic theorem, for $\mu$-a.e. $x\in X,$

$$\lim\limits_{n\to\infty} \frac{S_n(h,x)}{n}=\int\limits_X h d\mu$$

where the integral can be equal to infinity.

We write
 
$$ \frac{S_n(h,x)}{n}= \frac{S_n(f,x)}{n}+
\frac{g(T^n x)-g(x)}{n} $$

Since $g(x)$ is finite almost everywhere, we can choose a set $Y\subset X,$ such 
that $\mu(Y)>0$ and for every $y\in Y,$ $|g(y)|<M$ for some constant $M.$ Then by 
ergodicity of $T,$ there exists a subsequence $n_k,$ such that $T^{n_k}x\in Y,$ 
and therefore, by Birkhoff ergodic theorem for $\mu$- almost all $x\in X$ we have

$$\lim\limits_{k\to\infty} \frac{S_{n_k}(h,x)}{n_k}= \lim\limits_{k\to\infty}  
\frac{S_{n_k}(f,x)}{n_k}+ \frac{g(T^{n_k} x)-g(x)}{n_k}=0$$

which implies  $\int\limits_{X}h d\mu=0.$ 
\end{proof}


\begin{proof}[Proof of Theorem \ref{asexits}]
For any $\eps>0$ choose an $\alpha\in\Lambda(\eps),$ 
so that the map $\thtae^2$ is ergodic. Let $x\in I_\eps$ and $f(x)=n_{\eps}^1(x)-n_{\eps}^1(\thtae x).$  

Then the Birkhoff sums for $\thtae^2$ and  $f(x)$ are 

$$S_m(f,x)=f(x)+f(\thtae^2 x)+ f(\thtae^4 x)+ \ldots + f(\thtae^{2m} x)= $$
	$$=n_{\eps}^1(x)-n_{\eps}^2(x)+\ldots -n_{\eps}^{2m} (x)$$

By Proposition \ref{abs}, for Lebesgue almost every $x\in I_\eps,$ either there exists $m_0\in\N,$ such that
$S_{m_0}(f,x)\leq 0$ (and therefore $\qqe(x,\alpha)<\infty$)  or $f(x)$ is a co-boundary. 
But in the second case,  $S_m(f,x)=g(x)-g(\thtae^{2m+2} x)$ for a measurable $g(x).$
Either $g(x)< esssup\ {g(x)},$ or $g(x)=esssup\ {g(x)}$ on a positive Lebesgue measure set.
In either case,  the ergodicity of $\thtae^2,$  implies that for Lebesgue a.e. $x,$ there exists $m_0\in\N,$ 
such that $S_{m_0}(f,x)\leq 0$ and so  $\qqe(x,\alpha)<\infty.$

Now let $x\in S^1 \setminus I_\eps.$ Since $\alpha\notin\Q,$ there exists $n_0>0,$ such that
$R_{\alpha}^{-n_0}x\in I_\eps.$ Then for Lebesgue a.e. $x\in S^1$

$$
    \qqe(x,\alpha)\leq\qqe(\tht_{\alpha}^{-n_0}x,\alpha)<\infty
$$

\end{proof}

\section{Limiting distributions}
\label{distribution}

We now prove theorems \ref{thm1} and \ref{thm2}.

\subsection{Notations and the formulation of the main limiting distribution result}

Let $F_\eps^{(n)}(t_1,\ldots t_n)=\Prob\{\eps m_{\eps}^1>t_1,\
        \eps m_{\eps}^2>t_2,\ \ldots\ \eps m_{\eps}^n>t_n\}$ be the joint
        distribution function of the vector $\eps \bar{m}_{\eps}^n=(\eps m_\eps^1,\eps m_\eps^2,\ldots\eps
        m_\eps^n)^T.$

It is also convenient to introduce  
$$\nn_\eps(x,\alpha,T)=\#\{k\in\Z\cap (0,\eps^{-1}T] : k\alpha+x\subset I_\eps +\Z\},$$

the number of times the particle hits vertical walls during the time $\eps^{-1} T.$

Note that
\beq
	\Prob\{\eps m_\eps^k(x,\alpha)>t_k\ , k=1,\ldots, n\}=\Prob\{\nn_\eps(x,\alpha,t_k)\leq k-1\ , k=1,\ldots, n\}.
	\label{hittovisit}
\eeq

Let $\chi_I$ denote the characteristic function of the interval $I\subset\R$
and $\psi_{T}(x,y)=\chi_{(0,1]}(x)\chi_{[-T/2,T/2]}(y)$ be the characteristic function of a 
corresponding rectangle.
   
Then 
$$
        \nn_\eps(x,\alpha,T)=\sum\limits_{m=1}^{[\eps^{-1} T]}\sum\limits_{n\in\Z}\chi_{I_\eps}
        		\left( \alpha m+n+x\right) =$$

$$  =      \sum\limits_{(m,n)\in\Z^2}\chi_{(0,1]}\left(\frac{m}{\eps^{-1}T}\right)
            \chi_{[-T/2,T/2]}((\eps^{-1}T(\alpha m + n+x)) =
$$

$$
       = \sum\limits_{(m,n)\in\Z^2}\psi_{T}\left((m,n+x)\ua\dnet\right)
$$

Therefore, 

\beq
	\nn_\eps(x,\alpha,T)=\#\left\{ {(m,n)\in\Z^2} : (m,n+x)\ua\de \in\rect(T)  \right\},
\label{tolattice}
\eeq
where $\rect(T)=(0,T]\times [-1/2,1/2].$

Let $ASL(2,\R)=SL(2,\R)\ltimes \R^2$ be the semidirect product group with multiplication law
$$(M, {\bf v})(M',{\bf v}')=(M M', {\bf v} M'+{\bf v}').$$

The action of an element $(M, {\bf v})$ of this group on $\R^2$ is defined by 

\beq
	{\bf w} \mapsto {\bf w} M + {\bf v}
	\label{action}
\eeq

Each affine lattice of covolume one in $\R^2$ can then be represented as $\Z^2 g$ for some $g\in ASL(2,\R),$
and the space of affine lattices is represented by $X=ASL(2,\Z)\backslash ASL(2,\R),$ where
$ASL(2,\Z)=SL(2,\Z)\ltimes \Z^2.$ Denote by $\nu$ the Haar probability measure on $X.$

\begin{theorem}
	As $\eps\to 0,$ the limit of (\ref{hittovisit}) exists and is equal to 
	\begin{equation}
		F^{(n)}(t_1,\ldots t_n)=\nu(\{g \in X : \#\{\Z^2 g\cap\rect(t_k)\}\leq k-1 \ (k=1,\ldots, n)\}),
	\label{explicitlimit}
	\end{equation}
	which is a $C^1$ function $\R^n_{\geq 0}\to [0,1].$
	\label{limit}
\end{theorem}

We define the associated limiting probability density $\phi^{(n)}(t_1,\ldots, t_n)$ by

$$F^{(n)}(t_1,\ldots , t_n)=\int\limits_{t_1}^\infty\ldots \int\limits_{t_n}^\infty \phi^{(n)}(t_1,\ldots, t_n) dt_1,\ldots , dt_n$$

\subsection{The reduction of Theorem \ref{thm1} to Theorem \ref{limit}}

Because of the relation (\ref{relation}), Theorem \ref{limit} implies the convergence in 
distribution for the sequences $\{\eps n_\eps^k\}$ and $\{\eps \xi_\eps^k\}$ (part $(1)$ of Theorem \ref{thm1}).

Indeed,  let $k\in\N$ and $I_1,\ldots, I_k$ be a collection of $k$ intervals on the real line.
Let $I=I_1\times\ldots\times I_k\subset \R^k.$ Then 

$$ \lim\limits_{\eps\to 0} \Prob\{\eps n_\eps^1\in I_1,\ldots, \eps n_\eps^k\in I_k\}=
		\lim\limits_{\eps\to 0}\Prob\{\eps\bar{m}_\eps^k\in  {\bf B} I\}=$$
	$$		=\int\limits_{ {\bf B} I} \phi^{(k)}(t_1,\ldots, t_k)dt_1\ldots dt_k$$
	
and

$$ \lim\limits_{\eps\to 0} \Prob\{\eps\xi_\eps^1\in I_1,\ldots, \eps\xi_\eps^k\in I_k\}=
		\lim\limits_{\eps\to 0}\Prob\{\eps\bar{m}_\eps^k\in  {\bf B}{\bf A}^{-1} I\}=$$
	$$		=\int\limits_{ {\bf B}{\bf A}^{-1} I} \phi^{(k)}(t_1,\ldots, t_k)dt_1\ldots dt_k$$

The convergence for the random variable $\qqexa$ also follows from Theorem \ref{limit}.

Indeed, for any $k\geq 1$ let $\chi_{A_k}$ be the characteristic function of the set

$$  \Delta_k=\{(y_1,\ldots y_k)\in\R^k\ : y_1>0, \ldots, y_{k-1} >0, y_k<0\}. $$

Then for every $\eps>0$ we have
$$    \qqexa=\min\{j\in \Z_+:\ \xi_{\eps}^j\leq 0\} -1$$

Therefore,
$$
        \Prob\{\qqexa=k\}=\Prob\{\eps\xi_\eps^1>0,\ldots \eps\xi_\eps^{k-1}>0,\
           \eps \xi_\eps^k\leq0\}=$$
            
$$       =\Prob \{\eps \bar{m}_\eps^k \in {\bf B}{\bf A}^{-1} \Delta_k\} = \int\limits_{ {\bf B}{\bf A}^{-1}\Delta_k }dF_\eps^{(k)}, $$

and by Theorem \ref{limit} and the Helly-Bray Theorem (\cite{L}, p.183), there exists the limit

\beq
	G(k)=\lim\limits_{\eps\to 0}\Prob\{\qqexa=k\}=\int\limits_{ {\bf B}{\bf A}^{-1}\Delta_k } \phi^{(k)}(t_1,\ldots, t_k)dt_1\ldots dt_k
\label{limsum}
\eeq


Notice that the representation (\ref{limsum}) implies that $\sum\limits_{k=1}^{\infty} G(k)\leq 1$

\begin{proposition}

\begin{equation}
	\sum\limits_{k=1}^{\infty} G(k)=1
	\label{sumone}
\end{equation}
\end{proposition}

\begin{proof}
Let $\{\eta^k\}$ be the limiting process for the sequence $\{\eps n_\eps^k\}.$ From the explicit 
description (\ref{explicitlimit}) of the limiting distribution in Theorem \ref{limit}, we have the following 
description of the process $\{\eta^k\}.$

\begin{figure} [h]
	\centering
	\includegraphics [width=5in] {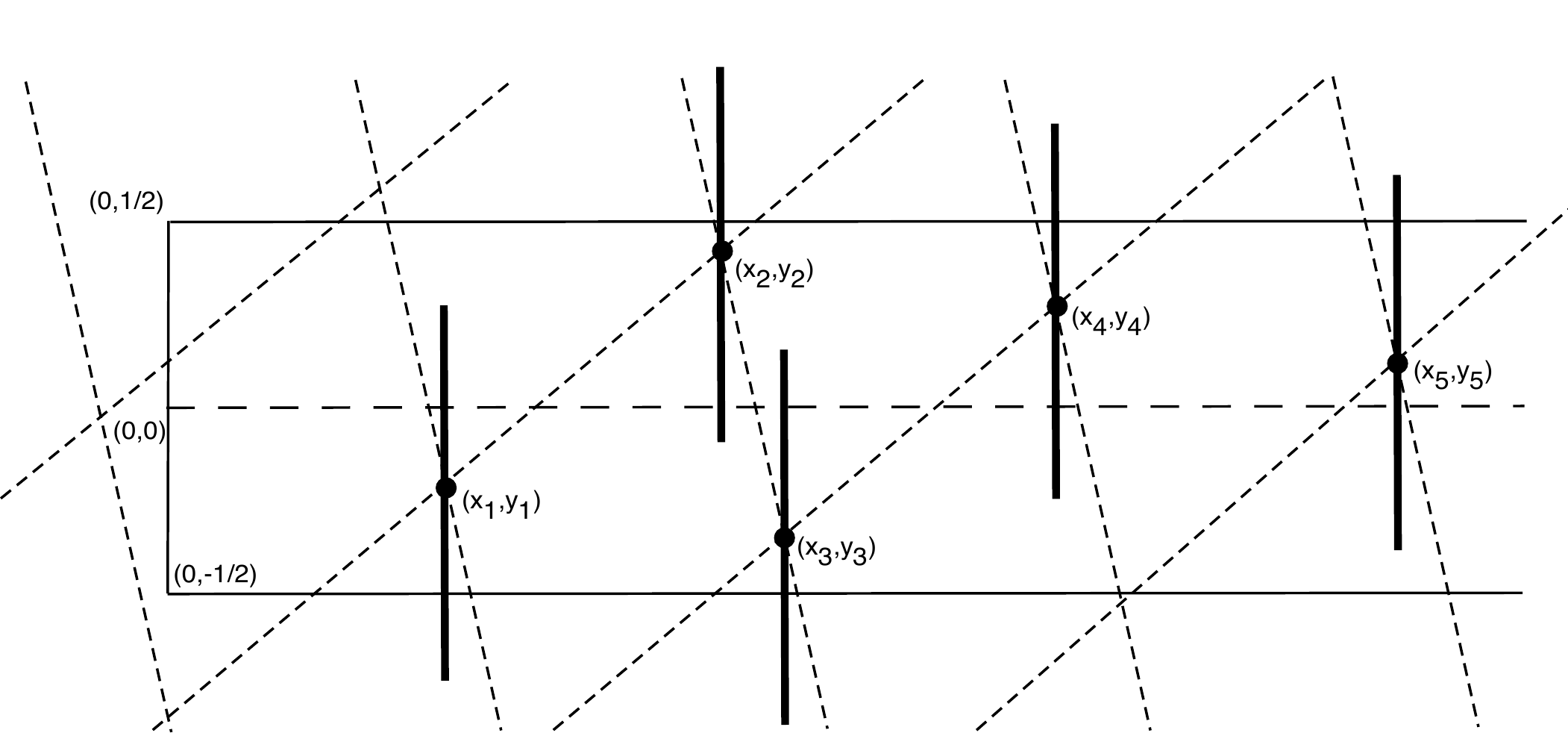}
	\caption{The horizontal ray through $(0,0)$ generates the sequence $\{-y_k\}$ as an orbit of an interval exchange map}
	\label{Fig2}
\end{figure}

Let $g \in X$ be an affine lattice which has no points either with the same horizontal coordinates, or 
on the boundary of the semi-infinite tube $\rect_\infty = [0,+\infty)\times [-1/2,1/2].$ The set of such lattices 
has full  Haar measure in $X.$ Let us enumerate points of $g$ which lie in $\rect_\infty$ according
to their horizontal coordinates: if the coordinates of the $k$'th lattice point of $g$ in $\rect_\infty$ are
$(x_k,y_k)=(x_k(g), y_k(g))$  $(k=1,2,\ldots),$  then  $x_k<x_{k+1}$ for any $k=1,2,\ldots$ Notice,
that $\nu$-almost every lattice $g$ has infinitely many points in $\rect_\infty.$


The sequence of random variables $\eps n_\eps^k=\eps n_\eps^k(x,\alpha)$ w.r.t. the probability
measure $\Prob$ on $\T^2$ converges in distribution to the sequence $\eta^1=\eta^1(g)=x^1(g),$ and
$\eta^k=\eta^k(g)=x_k(g)-x_{k-1}(g)$ for $k\geq 2$ w.r.t Haar measure $\nu$ on $X.$ 

Therefore, in order to prove (\ref{sumone}), it is enough to show that for $\nu$-almost every affine lattice
 $g\in X,$  there exists an even $k>0,$
such that 
\beq
	\eta^1-\eta^2+\eta^3- \ldots - \eta^k=x_1-(x_2-x_1)+(x_3-x_2)- \ldots -(x_k-x_{k-1}) \leq 0
	\label{limit_alt1}
\eeq

We will now show that the sequence $y_k(g)$ is an orbit of a certain map of an interval into itself, 
reduce (\ref{limit_alt1}) to a Birkhoff sum over this map and treat it in the way as in Section $2.$

First, we describe the map. Consider set $\mathcal{I}\subset\R^2$ of vertical segments of unit length centered 
at every lattice point of $g.$ We identify each segment in $\mathcal{I}$ with the base $I=[-1/2,1/2]$ of the tube 
$\rect_\infty$ by parallel translation. Let $\pi:\mathcal{I}\to I$ be the projection, which sends a point on some interval through a lattice point to the corresponding point in $I.$


Consider a unit speed flow in the positive horizontal direction on $\R^2.$ Its first return map to $\mathcal{I}$ is a well-defined map $\hat{T}=\hat{T}(g)$ of $\mathcal{I}$ into itself. We define the corresponding invertible map $T:I\to I,$ so that $\pi\circ\hat{T}=T\circ\pi.$ It is easy to see, that the map $T$ is an exchange of three intervals. For $\nu$-almost every lattice $g$ it has combinatorial type $(3\ 2\ 1).$

For every $y\in I,$ we let $\psi(y)$ to be the Euclidean distance between $\hat{y}\in\pi^{-1}(y)$ and its image under 
$\hat{T}.$  Clearly, this does not depend on the choice of $\hat{y}\in\pi^{-1}(y).$

Notice that the sequence $\{y_k\}$ of the vertical coordinates of the lattice points of $g$ in $\rect_\infty$ is related to the map $T$ described above: for $k\in\N$,  
$y_k=-T^{k-1}(-y_1)$ 
(see Figure 2). Also for $k\in\N,$ we have $\psi(-y_k)=x_{k+1}-x_k.$ Let $-y_0=T^{-1}(-y_1).$ Then the sum in (\ref{limit_alt1}) has the form (recall, $k$ is even)

\beq
	\begin{array}{c}
    		x_1-\psi(-y_1)+\psi(-y_2)-\ldots-\psi(-y_{k-1})\leq \\ \leq \psi(-y_0) - \psi(-T(-y_0))+\psi(T^2(-y_0))-\ldots-\psi(T^{k-1}(-y_0))
    \end{array}
\eeq


Therefore similarly to Section 2, the alternating sum (\ref{limit_alt1}) is reduced to a Birkhoff sum for the function $f(y)=\psi(-y)-\psi(-T(-y))$ and the map $T^2.$

 Let the lengths of the interval exchange map $T$ be equal to $(\lambda_1, \lambda_2, 1-\lambda_1-\lambda_2).$ Denote the 
simplex of possible $\lambda_i$'s by $$\Lambda=\{(\lambda_1,\lambda_2)\ |\ \lambda_1>0,\ \lambda_2>0,\ \lambda_1+\lambda_2<1\}\subset\R^2,$$ and the corresponding interval exchange map of combinatorial type $(3,2,1)$ by $T_{\lambda_1,\lambda_2}.$ 
The following theorem was first proved by Katok and Stepin in \cite{KS}.

\begin{theorem}
	For Lebesgue almost every pair $(\lambda_1,\lambda_2)\in\Lambda,$ the map $T_{\lambda_1,\lambda_2}$ of the interval 
	$I$ onto itself is weakly-mixing. 
\label{KS}
\end{theorem}

Similarly to the proof of Theorem \ref{asexits}, Theorem \ref{KS}  and Proposition \ref{abs} imply that there exists 
a full Lebesgue measure subset  $\Lambda_1\subset\Lambda,$ such that for every $(\lambda_1,\lambda_2)\in\Lambda_1,$
there exists  a full Lebesgue measure subset $I'=I'(\lambda_1,\lambda_2)\subset I,$ such that for every $y\in I'$ there 
exists $k>0,$ such that  $$\psi(-y_0) - \psi(-T_{\lambda_1,\lambda_2}(-y_0))+\psi(T_{\lambda_1,\lambda_2}^2(-y_0))-\ldots-\psi(T_{\lambda_1,\lambda_2}^{k-1}(-y_0))\leq 0.$$


Let $\tilde{X}\subset X$ be the set of lattices, for which the construction above gives an interval exchange transformation of combinatorial type
$(3\ 2\ 1).$ Then $\tilde{X}$ is open and $\nu(\tilde{X})=1.$ Notice that for any $g\in\tilde{X},$ the map $\mathcal{X}: g\mapsto (\lambda_1,\lambda_2, y_0)$ 
is differentiable and its differential is surjective. Therefore, the preimage of any Lebesgue measure zero set under $\mathcal{X}$ 
has Haar measure zero in $X.$ Therefore, the set of lattices $g\in X,$ such that $\mathcal{X}(g)\in\{(\lambda_1,\lambda_2, y_0)\ |\ (\lambda_1,\lambda_2)\in \Lambda_1,\ y_0\in  I'(\lambda_1,
\lambda_2)\}$ has full Haar measure in $X$ and so (\ref{sumone}) is proved.
\end{proof}

\begin{remark} The condition (\ref{sumone}) is equivalent to the tightness of the family of distributions $\{Q_\eps\}$ as $\eps\to 0.$ Namely,
for any $\delta>0$ there exists $N=N(\delta)$ and $\eps_1=\eps_1(\delta),$ such that for $\eps<\eps_1$
\beq
	1-\delta \leq \sum\limits_{k=1}^N\Prob\{\qqexa=k\} \leq 1
\label{tight}
\eeq 
\end{remark}

\subsection{Continuous case}

\begin{proposition}
\label{tail}
	For any $s>0$ and $\delta>0,$ there exists $\eps_0>0$ and $k\in \N,$ such that
	\beq
		\Prob\{(\eps m_{\eps}^{k}\leq s)\}< \delta
	\eeq  
	for all $\eps<\eps_0.$
\end{proposition}

\begin{proof}
	We have 
$$\Prob\{\eps m_\eps^k\leq s\}=\Prob\{\nn_\eps(x,\alpha,s)\geq k\}$$

The limit, as $\eps\to 0,$ exists and, in view of \cite{Ma2} (p.1131, first equation), is bounded by
$$\leq C_s k^{-3}$$ for some constant $C_s.$
\end{proof}

We now prove part $(1)$ of Theorem \ref{thm2}. 

For any $N\in\N$ and intervals $I_1,\ldots, I_N\subset \R,$
$$\Prob\{\eps\xi_\eps(s_1)\in I_1, \ldots, \eps\xi_\eps(s_N)\in I_N\}=$$
	$$=\sum\limits_{\bar{k}\in\Z_{\geq 0}^N}^{\infty}
		\Prob\{ \eps\xi_\eps(s_j)\in I_j,\ \eps m_\eps^{k_j}\leq s_j < \eps m_\eps^{k_{j+1}}\ (j=1,\ldots, N)\}$$

Notice that
$$\eps\xi_\eps(s)=\left\{
	\begin{array}{ll} 
			s & {\rm if}\ 0\leq s < \eps m_\eps^1,\\
			\eps \xi_\eps^k+(-1)^k (s-\eps m_\eps^k) & {\rm if}\ \eps m_\eps^k\leq s < \eps m_\eps^{k+1},
	\end{array}
	\right.$$
	
and 
$$\xi_\eps^k=\sum\limits_{i=1}^{k}(-1)^{i-1}(k-i+1)m_\eps^i.$$
	
Therefore, by Theorem \ref{limit}, for every {\it fixed} ${\bf {k}}\in\Z_{\geq 0}^N$ we have

$$\lim\limits_{\eps\to 0} \Prob\{ \eps\xi_\eps(s_j)\in I_j,\ \eps m_\eps^{k_j}\leq s_j < \eps m_\eps^{k_{j+1}}\ (j=1,\ldots, N)\}=$$
$$\int_{B_{\bf{k}}} \phi^{(k+1)}(t_1,\ldots,t_{k+1})dt_1\ldots dt_{k+1},$$

with $k=\max({\bf {k}}),$ and the range of integration restricted to the set 

\beq
	\begin{array}{lr}
		B_{\bf{k}}=\big\{ (t_1,\ldots, t_{k+1})\ :\qquad  t_{k_j} \leq s_j< t_{k_j+1},  \\ 
	       \sum\limits_{i=1}^{k}(-1)^{i-1}(k-i+1) t_i +(-1)^{k_j}(s_j-t_{k_j})\in A_j \big\}
	\end{array}
\eeq

Futhermore, 
$$\sum\limits_{\begin{array}{c} {\bf {k}}\in\Z_{\geq 0}^N\\ max({\bf {k}})\geq R\end{array}}^{\infty}
		\Prob\{ \eps\xi_\eps(s_j)\in I_j,\ \eps m_\eps^{k_j}\leq s_j < \eps m_\eps^{k_{j+1}}\ (j=1,\ldots, N)\}\leq $$
		
$$\sum\limits_{\begin{array}{c} {\bf {k}}\in\Z_{\geq 0}^N \\ k_1\geq R\end{array}}^{\infty}
		\Prob\{ \eps\xi_\eps(s_j)\in I_j,\ \eps m_\eps^{k_j}\leq s_j < \eps m_\eps^{k_{j+1}}\ (j=1,\ldots, N)\}\leq$$
$$\leq\Prob\{\eps m_\eps^R\leq s_1\}.$$

Part (1) of Theorem \ref{thm2} now follows from Proposition \ref{tail}.

For the part $(2)$ of Theorem \ref{thm2} we have,

\beq
	\Prob\{\eps T_\eps\leq s\}=\sum\limits_{k\in\N} \Prob\{\eps T_\eps\leq s,  \qqexa=k \}
	\label{limittime}
\eeq
	
Notice that if $\qqexa=k$, then the time which the particle spends in the tube is equal to 

$$
T_\eps=T_\eps(x,\alpha)= 2 \sqrt{1+\alpha^2} (n_\eps^1+n_\eps^3+\ldots + n_\eps^{k}),
$$	

 and so,
 
\beq
	\begin{array}{cr}
		\Prob\{\eps T_\eps\leq s,  \qqexa=k \}=\\ 		
		 =\Prob\{2 \eps\sqrt{1+\alpha^2} (n_\eps^1+n_\eps^3+\ldots + n_\eps^{k})<s,\  
			 \qqexa=k\}
	\end{array}
\label{cont_time_cond}
\eeq

By Theorem $8,$ for any $s>0$ there exists joint limiting distribution of
$$\Prob\{\alpha < s, \eps m_\eps^k(x,\alpha)>t_k\ (k=1,\ldots, n)\},$$
 as $\eps\to 0,$ and therefore, of (\ref{cont_time_cond}) as well.

On the other hand,
 $$\Prob\{\eps T_\eps\leq s,  \qqexa\geq k \} \leq \Prob\{\eps m^k_\eps\leq s\},$$

and so, Proposition \ref{tail} and the convergence of (\ref{cont_time_cond}) imply the existence of the limit
$$H(s)=\lim\limits_{\eps\to 0} \Prob\{\eps T_\eps\leq s\}$$

Also, since $$\Prob\{\eps T_\eps\leq s\}\geq \sum\limits_{k=1}^N \Prob\{\eps T_\eps\leq s,  \qqexa=k \},$$

the tightness (\ref{tight}) implies that $H(s)\to 1$ as $s\to\infty.$

This finishes the proof of part $(2)$ of Theorem \ref{thm2}.

\subsection{The proof of Theorem \ref{limit}}

By (\ref{hittovisit}) it is enough to show that for any $n\in\N$ and any $n$-tuples $(t_1,\ldots, t_n)\in\R^n_{>0},$
${\bf k}=(k_1,\ldots, k_n)\in\Z^n_{\geq 0}$ there exists the limit 
\beq
	\begin{array}{c}
	G^{(n)}(t_1,\ldots, t_n) = \lim\limits_{\eps\to 0} \Prob\{\nn_\eps(x,\alpha,t_j)= k_j,\ (j=1,\ldots, n)\}=  \\
	=\nu(\{g\in X : \#\{\Z^2 g\cap\rect(t_j)\}=k_j \ (j=1,\ldots, n)\}) 
	\end{array}
\label{limit2}
\eeq

and that $G^{(n)}(t_1,\ldots, t_n)$ is a $C^1$-function of $(t_1,\ldots, t_n).$

For $n=1$ the convergence in (\ref{limit2}) was first proved by Mazel and Sinai (\cite{MaSi}). It was later 
reproved and generalized by the third author (\cite{Ma1}, \cite{Ma2}) using different methods. The proof  of the 
convergence in (\ref{limit2}) follows the one in \cite{Ma1}. The proof of the regularity of the limiting function 
is similar to the one in \cite{MS}.

We reduce the convergence in (\ref{limit2}) to an equidistribution result for the geodesic flow on $X.$ 

Recall, that the action of the geodesic flow  on X is given by right action of a one-parameter subgroup of $X:$
$$
    \Phi^t=\left(\left(\begin{array}{cc} e^{-t/2} & 0\\ 0 & e^{t/2}\end{array}\right),(0,0)\right).
$$ 

The unstable horocycle of the flow $\Phi^t$ on $X$ is then parametrized by the subgroup 
$H=\{n_{-}(x,\alpha)\}_{(x,\alpha)\in\T^2}:$

$$
    n_{-}(x,\alpha)=\left(
    \left(\begin{array}{cc}
       1 & \alpha\\
       0 & 1
    \end{array}\right),
    (0,x)\right).
$$

For $g\in X$ let $F_T(g)$ be equal to the number of lattice points of $\Z^2 g$ in 
the rectangle $\rect(T).$

Then by (\ref{tolattice}) $$\nn_\eps(x,\alpha,T)=F_T(n_{-}(x,\alpha)\Phi^t)$$
with $t=-2\ln (\eps).$

\begin{theorem} \cite{Ma1} For any bounded $f: ASL(2,\Z) \backslash ASL(2,\R) \to\R,$ such that the
discontinuities of $f$ are contained in a set of $\nu$-measure zero and any
Borel probability measure $\Prob,$ absolutely continuous with respect to
Lebesgue measure on $[0,1)\times [0,1)$
$$
    \lim\limits_{t\to\infty}\int\limits_0^1\int\limits_0^1
    f(n_{-}(x,\alpha)\Phi^t)d\Prob(x,\alpha)= \int\limits_{ASL(2,\Z) \backslash ASL(2,\R)} f d\nu
$$
\label{equidistr}
\end{theorem}

Let 
$$D(g)=\left\{\begin{array}{ll} 
			1 & {\rm if}\ F_{t_j}(g)=k_j,\ (j=1,\ldots n),\\
			0 & {\rm otherwise}
		\end{array}\right.
$$

Then $D(g)$ satisfies the conditions of Theorem \ref{equidistr}. The
convergence in (\ref{limit2}) now follows from theorem \ref{equidistr} applied to the function $D(g).$

We now prove $C^1$ regularity of the limiting function
$G^{(n)}(t_1,\ldots,t_n).$ It is enough to consider the case when all
$t_j$ are different. We also assume that all $k_j>0.$ The case when
some $k_j=0$ is similar.


Let $X_1=SL(2,\Z)\backslash SL(2,\R)$ be the homogeneous space of
lattices of covolume one and let $\nu_1$ be the probability Haar
measure on $X_1.$ For a given ${\bf y}\in\R^2$ let
    $$X(\bfy)=\{g\in X\ :\ \bfy\in\Z^2 g\},$$
where the action of $X$ on $\R^2$ is given by the formula (\ref{action}).
    

There is a natural identification of the sets $X(\bfy)$ and $X_1$ through
$$X(\bfy)=\{(M,\bfy)\ :\ M\in X_1\}$$
Under this identification the probability Haar measure $\nu_1$ on
$X_1$ induces a probability Borel measure $\nu_\bfy$ on $X(\bfy).$

We will need the following two results.

\begin{proposition}(Siegel's formula, \cite{Sie})  Let $f\in L^1(\R^2),$ then
\beq
    \int\limits_{X_1}\sum\limits_{{\bf k}\in\Z^2\setminus{\bf 0}} f({\bf
    k}M) d\nu_1 (M)=\int\limits_{\R^2} f(x)dx
\eeq \label{Siegel}
\end{proposition}

\begin{proposition} (\cite{MS})
\label{Fubini2} Let $\mathcal{E}\subset X$ be any Borel set; then
$\bfy\to\nu_\bfy(\mathcal{E}\cap X(\bfy))$ is a measurable function
from $\R^2$ to $\R.$ If $U\subset\R^2$ is any Borel set such that
$\mathcal{E}\subset\cup_{\bfy\in U} X(\bfy),$ then
\beq
    \nu(\mathcal{E})\leq \int\limits_U \nu_\bfy(\mathcal{E}\cap
    X(\bfy))d\bfy
\label{Fubini} \eeq

Futhermore, if $\forall \bfy_1\neq\bfy_2 \in U\ :\ X(\bfy_1)\cap
X(\bfy_2)\cap\mathcal{E}=\emptyset,$ then equality holds in
(\ref{Fubini})

\end{proposition}

Notice that Propositions \ref{Siegel} and \ref{Fubini2} imply
that if there are two different indices $1\leq i,j\leq n,$ such that
$$
    \Delta_{ij} (h_i,h_j)= \{g\in X : |\Z^2 g\cap (\rect(t_i)\vartriangle\rect(t_i+h_i))|>0\}\cap
$$
$$
    \cap\{g\in X : |\Z^2
    g\cap(\rect(t_j)\vartriangle\rect(t_j+h_j))|>0\}\neq\emptyset,
$$

then
 $$
    \nu\{\Delta_{ij}(h_i,h_j)\}=\bar{o}(||{\bf h}||)\ {\rm as}\ ||h||\to 0
$$

Therefore,

\beq\begin{array}{cl}

G^{(n)}(t_1+h_1,\ldots, , t_n+h_n)-G^{(n)}(t_1,\ldots,t_n) = & \\

=\sum\limits_{j=1}^n G^{(n)}(t_1,t_2\ldots, t_{j-1},
t_j+h_j,t_{j+1},\ldots, t_n)-G^{(n)}(t_1,\ldots,t_n) +\bar{o}(||{\bf
h}||)= & \\

 =\sum\limits_{j=1}^{n} (\nu \{g\in X : |\Z^2 g\cap
\rect(t_j)|=k_j-1,\ |\Z^2 g\cap \rect(t_j+h_j)|=k_j, & \\

|\Z^2 g\cap \rect(t_i)|=k_i, i\neq j \}- \\

-\nu \{g\in X\ :\ |\Z^2 g\cap \rect(t_j)|=k_j,\ |\Z^2 g\cap
\rect(t_j+h_j)|=k_j+1\}, & \\

|\Z^2 g\cap \rect(t_i)|=k_i, i\neq j)+\bar{o}(||{\bf h}||) &

\end{array}
\label{long}
\eeq

Consider a single term in the expression above.

Let $$\mathcal{E}_j=\mathcal{E}_j(h_j)=\{g\in X : |\Z^2 g\cap
\rect(t_j)|=k_j,$$ $$|\Z^2 g\cap \rect(t_j+h_j)|=k_j+1, |\Z^2 g\cap
\rect(t_i)|=k_i, i\neq j\},$$

\medskip

and let $U=\rect(t_j+h_j)\setminus\rect(t_j).$ Then by the
proposition \ref{Fubini2},

$$\nu(\mathcal{E}_j)=\int\limits_U \nu_\bfy (\mathcal{E}_j \cap
X(\bfy))d\bfy=\int\limits_{t_j}^{t_j+h_j}\int\limits_{-1/2}^{1/2}
\nu_{(x,y)}(\mathcal{E}_j \cap X(x,y)) dx dy$$

Therefore, by proposition \ref{Siegel},

$$
	\lim\limits_{h_j\to 0}\frac1{h_j}\nu (\mathcal{E}_j(h_j))=
$$ 
$$
	= \int\limits_{-1/2}^{1/2}\nu_1(\{g\in X_1\ :\ |\Z^2 g\cap (\rect(t_i)-(t_j,y))|=k_i,\ (i=1,\ldots, n)\}) dy
$$

For every fixed $y\in [-1/2,1/2]$ continuity of the expression under
the integral sign with respect to $(t_1,\ldots,t_n)$ again follows
from Proposition \ref{Siegel}. It is clearly uniform in $y$ and
therefore the integral is continuous with respect to $(t_1,\ldots,
t_n).$ Each term in (\ref{long}) can be treated in a similiar way
and this proves $C^1$ regularity of the function $G^{(n)}(t_1,\ldots,
t_n)$ and finishes the proof of Theorem \ref{limit}.

\section{Proof of Theorem \ref{maintheorem}}

We now deduce (\ref{goal}) from part $(2)$ of Theorem \ref{thm1}.

Consider the unfolding of the tube to $\R^2$ obtained by the
reflections from the horizontal boundary of the tube. Let 
${\bf p_k}= (\xi_\eps^k,\zeta_\eps^k)$ be the position of the particle at the moment of
$k$'th reflection from the wall in this unfolding. 

\begin{figure} [h]
	\centering
	\includegraphics [width=4in] {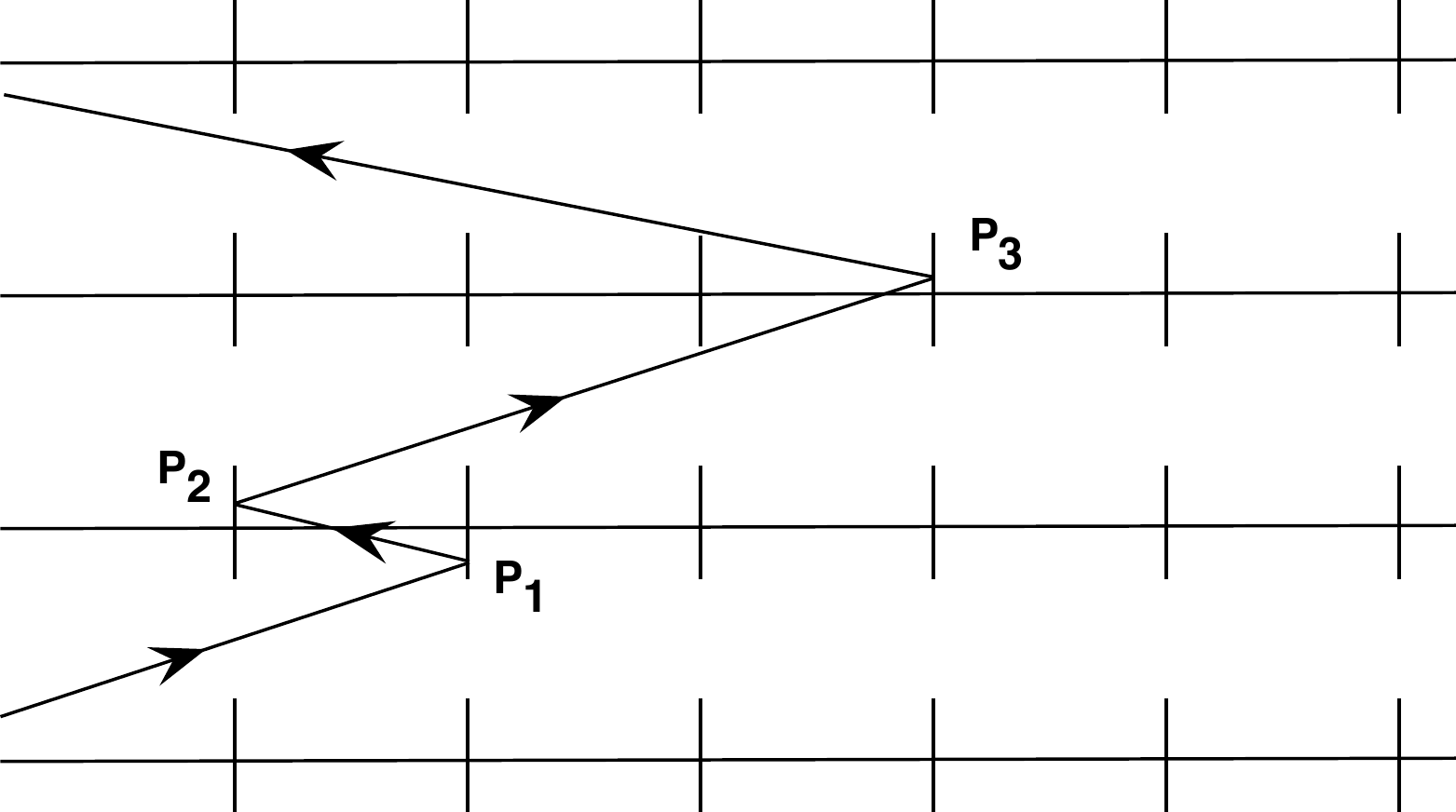}
	\caption{An unfolded trajectory. In this example, $\qq_\eps = 3,\, \lfloor \bar{\zeta} \rfloor =2$ and $n^1_\varepsilon = 2,\, n^2_\varepsilon = 1,\ n^3_\varepsilon =3,\, n^4_\varepsilon > 4$}
	\label{Fig3}
\end{figure}

Then
$$\xi_\eps^k=n_\eps^1-n_\eps^2+\ldots+ (-1)^{k+1}n_\eps^k$$
and
$$\zeta_\eps^k=y_{\rm in}+\alpha(n_\eps^1+n_\eps^2+\ldots+n_\eps^k)$$

At the moment of the exit from the tube, the
vertical coordinate of the particle is

\beq
  \bar{\zeta}=2(y_{\rm in}+\alpha n^1_\eps+\alpha n^3_\eps+\ldots+ \alpha n^{\qq_\eps}
_\eps)-y_{\rm in}
\label{ineq0}
\eeq
Let 
$$z=y_{\rm in}+\alpha n^1_\eps+\alpha n^3_\eps+\ldots+ \alpha n^{\qq_\eps}
_\eps
$$

and let 
$||\cdot||$ denote the distance to the nearest integer. 

Then
$$||y_{\rm in}+\alpha n_\eps^1||\leq \eps/2, ||\alpha n_\eps^i||\leq \eps\ {\rm for}\ i>1$$

Therefore,
\beq
    || z ||\leq \frac{\eps \qq_\eps}{2}
\label{ineq}
\eeq

Notice that $v_{\rm out}=-v_{\rm in},$ if both the number of reflections from
vertical walls and from horizontal walls is odd. The former is obviously
odd at the moment of exit. The number of reflections from the horizontal walls  
is equal to the integer part $\lfloor \bar{\zeta}\rfloor.$

If $z-\lfloor z\rfloor \leq 1/2,$ then by (\ref{ineq0}), $\lfloor \bar{\zeta}\rfloor$  is odd provided that $2||z||<y_{\rm in},$
and if $z-\lfloor z\rfloor > 1/2,$ then $\lfloor \bar{\zeta}\rfloor$  is odd provided that $1-2||z||>y_{\rm in}.$

 By (\ref{ineq}) this is the case, when
$$
	\eps \qq_\eps<\min\{y_{\rm in},1-y_{\rm in}\}
$$

By the assumption, the probability measure $\Prob$ on the initial conditions $(y_{\rm in},\alpha)$
is absolutely continuous with respect to the Lebesgue measure, therefore for any $k\in\N,$ 
$$		\Prob\{\qq_\eps=k,\ \eps k < \min\{y_{\rm in},1-y_{\rm in}\}\}  = $$
$$		\Prob\{\qq_\eps=k\}-\Prob\{\qq_\eps=k,\   \min\{y_{\rm in},1-y_{\rm in}\} \leq \eps k \}  \to G(k)\ {\rm as}\ \eps\to 0$$

Together with the tightness condition (\ref{tight}) this implies 
$$
    \Prob\{\eps \qq_\eps <\min\{y_{\rm in},1-y_{\rm in}\}\}\to 1\ {\rm as}\ \eps\to 0
$$
and so,
$$\Prob\{v_{\rm out}=-v_{\rm in}\} \to 1\ {\rm as}\ \eps\to 0.$$
Note that the existence of the limiting probability distribution for $\{\qq_\eps\}$ as $\eps\to 0$  also implies that for any
$\delta>0,$
$$\Prob\{|y_{\rm out}-y_{\rm in}|>\delta\}\to 0\ {\rm as}\ \eps\to 0.$$
Indeed, after each reflection from a vertical wall, the particle backtracks itself with an error at most $\eps,$
so at the moment of exit it backtracks the incoming trajectory with total error of at most $\eps\qq_\eps.$

This finishes the proof of Theorem \ref{maintheorem}.


\begin{thebibliography}{GGOR}

	\bibitem{B} Boshernitzan, M.,  {\it A condition for minimal interval exchange maps to be uniquely ergodic,} Duke Math. J. {\bf 52} (1985) pp. 723--752,

	\bibitem{BN} Boshernitzan, M., Nogueira, A., {\it Generalized functions of interval exchange maps,} Ergodic Theory and Dynamical Systems, {\bf 24 } (2004) pp. 697--705,

	\bibitem{E}   Eaton, J.E.,  {\it On spherically symmetric lenses,}  Trans. IRE Antennas Propag. {\bf 4} (1952) pp. 66--71,

	\bibitem{KS}  Katok, A., Stepin, A., {\it Approximations in Ergodic Theory,} Russian Math. Surveys, {\bf 22} (1967) n. 5, pp. 77--102,  
	
	\bibitem{L}   Loeve, M.,  {\it Probability Theory I,} Springer-Verlag, Berlin-Heidelberg-New York,
    1977
	
	\bibitem{Ma1} Marklof, J. Distribution modulo one and Ratner's theorem, {\it Equidistribution in Number Theory, An Introduction,} eds.
    A. Granville and Z. Rudnik, Springer 2007, pp. 217-244,

	\bibitem{Ma2} Marklof, J. {\it The $n$-point correletions between values of a linear form}, Ergodic Theory and Dynamical Systems, {\bf 20} (2000), pp. 1127--1172,

	\bibitem{MS}   Marklof, J., Str\"{o}mbergsson, A., {\it The distribution of free path lenghts in the periodic Lorentz gaz
    and related lattice point problems}, to appear in the Annals of Mathematics,
	
	\bibitem{MaSi} Mazel, A.E., Sinai Ya.G., {\it A limiting distribution connected with fractional parts of linear forms}, in: Ideas and Methods in Mathematical Analysis, Stochastics and Applications, S. Albeverio {\it et al.} (eds.), Cambridge Univ.
Press, Cambridge, 1992, pp. 220--229,

	\bibitem{PG} Plakhov, A.,  Gouveia, P.,  {\it Problems of maximal mean resistance on the plane}, Nonlinearity {\bf 20} (2007), pp. 2271-2287,
	
         \bibitem{Sie}  Siegel, C.L., {\it Lectures on the Geometry of  Numbers}, Springer-Verlag, Berlin-Heidelberg-New York, 1989
    
         \bibitem{T} Tyc, T., Leonhardt, U., {\it Transmutation of singularities in optical instruments,} New J. Physics {\bf 10} (2008) 115038 (8pp)

\end{thebibliography}
\end{document}